# Image and Transfer Functions

**Dr. Emil Schwab, Professor**
"Tibiscus" University of Timisoara

**ABSTRACT.** We describe three transfer functors $P, P', P''$ of an inverse exact category which arise from three transfer functions. We concentrate on some of the basic results which emerge from the theory of projections in inverse exact categories.

## 1. Baer*-categories and exact inverse categories

Inverse categories have been considered by a number of authors: J.Kastl [Kas79], M.Grandis [Gra76], H.J.Hoehnke [Hoe88], M.V.Lawson [Law91] and by the author of this paper ([Sch81]). A category C is said to be inverse if for each morphism f there exists a unique morphism g such that

$$fgf = f \quad \text{and} \quad gfg = g.$$

The morphism g is called the inverse of f and we shall denote it by $f^{-1}$. The inverse is then nothing but an involution of the inverse category called the canonical involution. An involution of a category is a contravariant endofunctor identical on objects and involutory on morphisms. In a category with an involution $*$, if for a morphism f there exists a morphism g such that

$$fgf = f, \quad gfg = g, \quad (fg)^* = fg, \quad (gf)^* = gf$$

then it is unique and it is called the Moore-Penrose generalized inverse of f with respect to $*$. The generalized inverse of a morphism f is denoted by $f^{(-1)}$.

An idempotent morphism i such that i*=i is called projection. A category B with zero object and with involution ∗ is called a Baer*-category if for each morphism f there exists a (unique) projection $f'$ so that

$$\{g \in B \mid fg = 0\} = f'B.$$

A projection i of a Baer*-category is called closed if $(i')' = i$ (often, $(f')'$ is denoted by $f''$).

An exact category is a normal and conormal category with kernels and cokernels and with the property that every morphism f can be written as a composition f=pq where p is a monomorphism and q is an epimorphism.

The main result of this section is the following one:

**Theorem 1.1.** *Let* C *be a category in which any morphism has a Moore-Penrose generalized inverse with respect to an involution* ∗. *Then* C *is exact if and only if* C *is a Baer*-category with closed projections so that every projection i can be written as a composition i=pq where p is a monomorphism and q is an epimorphism.*

**Proof.** Assume that C is exact. Let f be a morphism of C and u=kerf. Then, the morphism $i=uu^{(-1)}$ is a projection such that fi=0. Now, if g is a morphism of C such that fg=0, then there exists a morphism h such that g=uh and therefore ig=g. It follows that {g∈C|fg=0}=iC, that is:

$$f' = \ker f (\ker f)^{(-1)}$$

for any morphism f of C. Consequently, C is a Baer*-category.
Since C is exact, any projection i has a mono-epi factorization: i=pq, and any monomorphism is the kernel of a morphism. It follows that p=kerh for some morphism h of C. Then hi=0, and for any morphism g of C such that hg=0, there exists a morphism k of C such that g=pk. Consequently,

$$g = pk = pqpk = i(pk) = ig.$$

This shows that

$$i = h'.$$

Now, it is straightforward to check that for any morphism h of a Baer*-category we have $h' = \{[h']'\}'$. Thus, it follows that any projection i of a Baer*-category is closed.

Conversely, if the category C with Moore-Penrose inverses is a Baer*-category which closed projections so that every projection can be written as a composition of a monomorhpism with an epimorphism, then for any morphism of C, the mono-epi factorizations

$$f' = p_1 q_1 \quad and \quad (f^*)' = p_2 q_2$$

implies:

$$p_1 = \ker f \quad and \quad q_2 = co\ker f.$$

Thus, C is a category with kernels and cokernels. To show that C is normal and conormal, let u be a monomorphism of C and let v be an epimorphism of C. It is straightforward to check that

$$u = \ker(u^*)' \quad and \quad v = co\ker v'.$$

So, C is normal (every monomorphism is the kernel of some morphism) and conormal (every epimorphism is the cokernel of some morphism).
Now, let f be a morphism of C and let

$$[(f^*)']' = pq$$

be a mono-epi factorization of the projection $[(f^*)']'$. Then the monomorphism p is an epimorphic image of f.
Consequently, C is an exact category.

The canonical involution * of an inverse category is defined by: $f^* = f^{-1}$. It is easy to see that $f^{-1}$ is the Moore-Penrose inverse of $f$ with respect to the canonical involution of an inverse category. So, an inverse category is exact if and only if it is a Baer*-category with closed projections in which every projection has a mono-epi factorization.

The category of partial bijections between sets is an exact inverse category. Another exact inverse category is the following one: the objects are all the sets with base points; the morphisms from (A,a) to (B,b) are the set functions f from A to B such that for any y∈B there exists x∈A so that f(x)=y and f(a)=b; the composition is the usual composition of maps. An inverse monoid S with zero adjoined is an inverse category with two objects. This inverse category is exact if and only if S is a group.

## 2. The image functor P

Transfer functions were introduced by Grandis [Gra77]. The orthodox expansion of a regular category with involution is constructed based on transfer functions. If R is a regular category (i.e. for any morphism f there is a morphism g such that fgf=f) with an involution $*$, and f is a morphism of R from A to B then the function $T(f): Hom_R(A,A) \to Hom_R(B,B)$ defined by

$$T(f)(h) = fhf^* \quad \text{for any } h \in Hom_R(A,A)$$

is called the transfer function for f.

Now, let C be an exact inverse category and let $*$ be the canonical involution on C. If f is a morphism of C from A to B we call the restriction on P(A) (the set of all projections from A to A) of the transfer function T(f), the image function for f. The functor P:C→Ens defined by:

$$A \in ObC \longrightarrow P(A) \quad ; \quad f \in Hom_C(A,B) \longrightarrow P(f) \in Hom_{Ens}(P(A), P(B)),$$

where P(f):P(A)→P(B) is the image function for f:

$$P(f)(i) = fif^* \quad (i \in P(A))$$

is a covariant functor.

**Theorem 2.1** *If* u∈Hom<sub>C</sub>(X,A) *is a monomorphism*, f∈Hom<sub>C</sub>(A,B) *is a morphism and*

$$P(f)(uu^*) = pp^*$$

*with* p:I→B *monomorphism, then p is the image of* fu (i.g. f(X)=I).

148

**Proof.** We have:

$$fu = fu(fu)*fu = fuu*f*fu = P(f)(uu*)fu = pp*fu = pq$$

where q=p*fu. Now, let fu=st a factorization of fu with s monomorphism. Then,

$$ss*pp*fu = ss*fu = ss*st = st$$

and therefore,

$$t = s*pp*fu.$$

It follows:

$$ss*pp* = ss*P(f)(uu*) = ss*fuu*f* = ss*pp*fuu*f* = stu*f* = fuu*f* = pp*$$

Consequently:

$$ss*p=p.$$

Thus, (I,p) is the smallest subobject of B which fu factors through. So, p is the image of fu.

Taking into account Theorem 2.1 we say that P is the image functor of the exact inverse category C. Some of the elementary properties of the image functor appear in the next theorem.

**Theorem 2.2** *Let* C *be an exact inverse category and* P *the image functor of* C. *The following properties are true:*
   (i)    P *preserve monomorphisms and epimorphisms*;
   (ii)   P(f)(0)=0 *and* P(f)(1)=ff* *for any morphism* f *of* C;
   (iii)  P(f)(f*f)=ff* *for any morphism* f *of* C;

**Proof.** (i) Let f be a monomorphism of C from A to B. If $P(f)(i) = P(f)(j)$ (where $i, j \in P(A)$), that is $fif^* = fjf^*$, then $if^* = jf^*$ since f is monomorphism. This result implies

$$(fi)^* = (fj)^* \Rightarrow fi = fj \Rightarrow i = j$$

and so, $P(f)$ is a monomrphism of Ens.

Now, let f be an epimorphism of C from A to B and $j \in P(B)$. Since $ff^* = 1_B$ it follows

$$(f^* jf)(f^* jf) = (f^* jf) \quad and \quad (f^* jf)^* = f^* jf.$$

Hence $f^* jf \in P(A)$ and so

$$P(f)(f^* jf) = ff^* jff^* = j.$$

Thus $P(f)$ is an epimorphism of Ens.

(ii) and (iii) are obviously.

An inverse semigroups in which every element is an idempotent are precisely the meet semilattices. If S is an inverse semigroup in which every element is an idempotent then a relation ≤ on S defined by

$$e \leq f \Leftrightarrow e = ef (= fe)$$

is a partial order on S and $e \cap f = ef$. Conversely, let (S,≤) be a meet semilattice. It is routine to check that S is a commutative semigroup with respect to the operation $\cap$, the greatest lower bound. Clearly, $e = e \cap e$ for each element e∈S. Thus (S, $\cap$) is an inverse semigroup in which every element is idempotent.

In an exact inverse category C for any object A, P(A) is an inverse semigroup in which every element is an idempotent. Thus P(A) is a meet semilattice. We have:

**Theorem 2.3.** *Let C be an exact inverse category and let f be a morphism of C from* A *to* B. *Then:*
   (i)   *P(f) is an inverse semigroups homomorphism;*

(ii) *P(f) has the order preserving property;*
(iii) *P(f)(i)≤ff\* for any i∈P(A)*
(iv) *P(f)(i)=ff\* if i≥f\*f*

**Proof.** (i).
$$P(f)(ij) = fijf^* = fif^* fjf^* = P(f)(i) \cdot P(f)(j) \quad \text{for any } i, j \in P(A).$$

(ii).
$$i \leq j \implies ij = i \implies P(f)(i) \cdot P(f)(j) = fif^* fjf^* = fijf^* = fif^* = P(f)(i) \implies$$
$$\implies P(f)(i) \leq P(f)(j).$$

(iii).
$$P(f)(i) \cdot ff^* = fif^* ff^* = fif^* = P(f)(i) \implies P(f)(i) \leq ff^*.$$

(iv).
$$i \geq f^*f \implies if^*f = f^*f \implies P(f)(i) \cdot ff^* = fif^* ff^* = ff^* ff^* = ff^* \implies$$
$$\implies ff^* \leq P(f)(i).$$

By (iii) it is now clear that if $i \geq f^*f$ then $P(f)(i) = ff^*$.

## 3. The inverse image functor $P'$

In the case of exact inverse categories, we have another "transfer functions". Let C be an exact inverse category and let f be a morphism of C from A to B. We call the function $P'(f): P(B) \to P(A)$ defined by

$$P'(f)(j) = (j'f)' \quad (j \in P(B))$$

the inverse image function for f. The functor $P': C \to Ens$ defined by

$$A \in ObC \longrightarrow P(A) \quad ; \quad f \in Hom_C(A, B) \longrightarrow P'(f) \in Hom_{Ens}(P(B), P(A))$$

is a contravariant functor.

**Theorem 3.1.** *If* $v \in Hom_C(A,B)$ *is a monomorphism,* $f \in Hom_C(A,B)$ *is a morphism and*

$$P'(f)(vv^*) = pp^*$$

*with* $p : \cdot \to A$ *monomrphism, then the diagram*

$$\begin{array}{ccc} \cdot & \xrightarrow{v^* fp} & \cdot \\ \downarrow p & & \downarrow v \\ A & \xrightarrow{f} & B \end{array}$$

*is a pullback.*

**Corollary 3.2.** *Let*

$$\begin{array}{ccc} A & \xrightarrow{f} & B \\ \uparrow u & & \uparrow v \\ X & & Y \end{array}$$

be a diagram of C with u and v monomorphisms. Then

(i)  $f(X) = Y$  *if and only if*   $P(f)(uu^*) = vv^*$
(ii) $f^{-1}(Y) = X$  *if and only if*   $P'(f)(vv^*) = uu^*$

*where f(X) denotes the image of the composition* $X \xrightarrow{u} A \xrightarrow{f} B$ *and* $f^{-1}(Y)$ *is the inverse image of Y in the sense of Mitchell [Mit65].*

The basic structural properties of $P'$ are similar with the properties (see Theorems 2.2. and 2.3.) of $P$ :

**Theorem 3.3.** *Let* C *be an exact inverse category and* $P'$ *the inverse image functor of* C. *The following properties are true:*
(i)   $P'(f)$ *is a monomorphism (epimorphism) if and only if f is an epimorphism (monomorphism)*
(ii)  $P'(f)(0) = f'$ *and* $P'(f)(1) = 1$ *for any morphism* f *of* C;
(iii) $P'(f)(ff^*) = 1$ *for any morphism* f *of* C;

**Theorem 3.4.** *Let C be an exact inverse category and let f be a morphism of C from* A *to* B. *Then:*
  (i)   $P'(f)$ *is an inverse semigroups homomorphism;*
  (ii)  $P'(f)$ *has the order preserving property;*
  (iii) $P'(f)(j) \geq f'$ *for any* $j \in P(B)$
  (iv)  $P'(f)(j) = 1$ *if* $j \geq ff^*$

The following properties are connection properties between $P$ and $P'$:

**Theorem 3.5.** *let C be an exact inverse category. Then*
  (i)   $P'(f) = P(f^*)$ *if and only if f is a monomorphism;*
  (ii)  $P(f) = P'(f^*)$ *if and only if f is an epimorphism;*
  (iii) $P(f)P'(f)P(f) = P(f)$ *and* $P'(f)P(f)P'(f) = P'(f)$ *for any morphism f of C.*

## 4. The transfer funtor $P''$

We now turn to another transfer function in the development of a new transfer functor $P''$. If C is an inverse exact category, then the functor $P'': C \to Ens$ defined by

$$A \in ObC \longrightarrow P(A) \quad ; \quad f \in Hom_C(A,B) \longrightarrow P''(f) \in Hom_{Ens}(P(B), P(A))$$

where

$$P''(f)(j) = (jf)''$$

is a contravariant functor.

Some of the elementary properties of $P''$ appear in the next theorems.

**Theorem 4.1.** *Let* C *be an exact inverse category and* $P''$ *the transfer functor defined above. The following properties are true:*
  (i) $P''(f)$ *is a monomorphism (epimorphism) if and only if f is an epimorphism (monomorphism)*

(ii)    $P''(f)(0) = 0$ and $P''(f)(1) = f''$ for any morphism f of C;
(iii)   $P''(f)((f*)') = 0$ for any morphism f of C;

**Theorem 4.2.** *Let C be an exact inverse category and let f be a morphism of C from* A *to* B. *Then:*
*(v)*    $P''(f)$ *is an inverse semigroups homomorphism;*
*(vi)*   $P''(f)$ *has the order preserving property;*
*(vii)*  $P''(f)(j) \leq f''$ *for any $j \in P(B)$*
*(viii)* $P''(f)(j) = 0$ *if $j \leq (f*)'$*

We note that between the transfer functions defined above they are interesting connections. For example, it is easy to establish that for any morphism $f : A \to B$ of the inverse exact category C, we have:

$$P''(f)(j) = (P'(f)(j'))' \qquad (\forall j \in P(B)).$$

This connection give rise equivalences between Theorems of Sections 3 and 4 (for example: Theorems 3.3 and 4.1., or 3.4. and 4.2.). We prove here some such equivalences:

(i) $P'(f)$ *is a monomorphism (epimorphism)* $\Leftrightarrow$ $P''(f)$ *is a monomrphism (epimorphism)*

($\Rightarrow$) If $P'(f)$ is a monomorphism then:

$P''(f)(j_1) = P''(f)(j_2) \Rightarrow (P'(f)(j'_1))' = (P'(f)(j'_2))' \Rightarrow P'(f)(j'_1) = P'(f)(j'_2) \Rightarrow$

$\Rightarrow j'_1 = j'_2 \Rightarrow j_1 = j_2$.

($\Leftarrow$) If $P''(f)$ is a monomorphism then:

$P'(f)(j_1) = P'(f)(j_2) \Rightarrow (P'(f)(j_1))' = (P'(f)(j_2))' \Rightarrow P''(j'_1) = P''(j'_2) \Rightarrow$

$\Rightarrow j'_1 = j'_2 \Rightarrow j_1 = j_2$.

(ii) $P'(f)(1) = 1 \Leftrightarrow P''(f)(0) = 0$;

154

($\Rightarrow$) $P'(f)(1) = 1$ implies:

$$P''(f)(0) = (P'(f)(0'))' = (P'(f)(1))' = 1' = 0.$$

($\Leftarrow$) $P''(f)(0) = 0$ implies:

$$(P'(f)(1))' = (P'(f)(0'))' = P''(f)(0) = 0 \Rightarrow P'(f)(1) = 1.$$

(iii) $P'(f)(0) = f' \Leftrightarrow P''(f)(1) = f''$.

($\Rightarrow$) $P'(f)(0) = f'$ implies:

$$P''(f)(1) = (P'(f)(1'))' = (P'(f)(0))' = f''.$$

($\Leftarrow$) $P''(f)(1) = f''$ implies:

$$(P'(f)(0))' = (P'(f)(1'))' = P''(f)(1) = f'' \Rightarrow P'(f)(0) = f'.$$